\documentstyle[12pt]{article}
\def \Z{\hbox{$Z\hskip -5.2pt Z$}}
\def \cc{\hbox{$c\hskip -4.6pt c$}}
\def\sZ{\hbox{$\sc Z\hskip -4.2pt Z$}}
\def \C{\hbox{$C\hskip -5pt \vrule height 6pt depth 0pt \hskip 6pt$}}
\def\qed{\hfill \hfill \ifhmode\unskip\nobreak\fi\ifmmode\ifinner
         \else\hskip5pt\fi\fi
 \hbox{\hskip5pt\vrule width4pt height6pt depth1.5pt\hskip 1 pt}}
\def\a{\alpha}
\def\e{\epsilon}

\def\d{\delta}
\def\D{\Delta}

\def\l{\lambda}
\def\L{\Lambda}

\def\si{\sigma}

\def\sc{\scriptstyle}
\def\ssc{\scriptscriptstyle}
\def\dis{\displaystyle}
\def\cl{\centerline}

\def\ol{\overline}
\def\ul{\underline}
\def\wt{\widetilde}
\def\wh{\widehat}
\def\rar{\rightarrow}
\def\Rar{\Rightarrow}
\def\lar{\leftarrow}
\def\Lar{\Leftarrow}

\def\Lra{\Leftrightarrow}
\def\bs{\backslash}
\def\hs{\hspace*}
\def\vs{\vspace*}
\def\rb{\raisebox}
\def\ra{\rangle}
\def\la{\langle}
\def\ni{\noindent}
\def\hi{\hangindent}
\def\ha{\hangafter}
\begin{document}
\par\
\par\
\par\ni
                 CLASSIFICATION OF INFINITE DIMENSIONAL WEIGHT
                 MODULES OVER THE LIE SUPERALGEBRA $sl(2/1)$
\par\
\par\ni
                 Yucai Su
\par\ni
                 Department of Applied Mathematics,
                 Shanghai Jiaotong University, China
\par\
\par\ni
                 ABSTRACT.
We give a complete classification of infinite dimensional indecomposable
weight modules over the Lie superalgebra $sl(2/1)$.
\par\
\par
\cl{             \S1. Introduction
}
\par
Among the basic-classical Lie superalgebras classified by Kac [3], the
lowest dimensional of these is the Lie superalgebra $B(0,1)$ or $osp(1,2)$,
while the lowest dimensional of these which has an isotropic odd simple
root is the Lie superalgebra $A(1,0)$ or $sl(2/1)$. The finite dimensional
indecomposable $B(0,1)$-modules are known to be simple by Kac [4].
Chmelev [1] classified finite dimensional indecomposable weight
$sl(2/1)$-modules and Leites [5] generalized the result to $sl(m/1)$.
We found some indecomposable generalized weight $sl(2/1)$-modules in [7]
(thanks are due to Germoni [2], who pointed out that the list in [7] is
incomplete). More generally, [2] 
obtained the classification of indecomposable (weight or generalized weight)
$sl(m/n)$-modules of singly atypical type and [2] proved that
indecomposable $sl(m/n)$-modules of other types are wild and are unable
to be classified.
\par
In [6], by using the technique employed in [8] by us in the classification
of indecomposable $sl(2)$-modules, we were able to give the classification
of infinite dimensional indecomposable $B(0,1)$-modules. In this paper, we
give the classification of infinite dimensional indecomposable weight
$sl(2/1)$-modules (but not for generalized weight modules) and obtain our
main result in Theorem 3.5.
\par
The problem of classifying infinite dimensional indecomposable modules of
the basic-classical Lie superalgebras has so far received little attention
in the literature. However, we believe that our classification of infinite
dimensional indecomposable $sl(2/1)$-modules will certainly help us to
better understand modules over (finite or infinite dimensional) Lie
superalgebras, just as we have seen in [8] that the classification of
infinite dimensional indecomposable $sl(2)$-modules has helped us to
understand modules over the Virasoro algebra. This is our motivation to
present the results here.
\par\
\par
\cl{              \S2. Notations and preliminary results
}
\par
First let us recall some basic concepts. Let $G$ denote the Lie superalgebra
$sl(2/1)$, consisting of $3\times3$ matrices of supertrace 0, which can be
defined as $G=G_{\ol 0}\oplus G_{\ol 1}$ with
$G_{\ol 0}=\la e_1,f_1,h_1,h_2\ra$, $G_{\ol1}=\la e_2,f_2,e_3,f_3\ra$,
where $G_{\ol0}$ and $G_{\ol1}$ are respectively the even and odd subspaces.
Let $(G_{\ol0})_{\pm},H,G_\pm,G_{\pm1}$ have their usual meanings. As a
basis for $G$, one can take
$$
\matrix{
e_1=E_{12},\hfill&e_2=E_{23},\hfill&e_3=E_{13},\hfill&h_1=E_{11}-E_{22},
\vs{2pt}\hfill\cr
f_1=E_{21},\hfill&f_2=E_{32},\hfill&f_3=E_{31},\hfill&h_2=E_{22}+E_{33},
\hfill\cr
}
\eqno(2.1)$$
where $E_{ij}$ is the $3\times3$ matrix with entry 1 at $(i,j)$ and 0
otherwise. Then $H=\la h_1,h_2\ra$. For convenience, set
$h_3=h_1+h_2$. Let $\a_1=(\e_1-\e_2)|_H,\a_2=(\e_2-\e_3)|_H$ be
the simple roots, where $\e_i(E_{jj})$ $=\d_{ij}$, then
$\{a_1,\a_2\}$ is a basis of $H^*$ and $D^+_1=\{\a_2,\a_1+\a_2\}$
is the set of positive roots of $G_1$, and
$S=\D^+_1\cup\{\a_1+2\a_2\}$ that of the universal enveloping
algebra $U(G_1)$. Define $\rho_0,\rho_1$ to be respectively half the sum of
all positive even,\,odd roots, $\rho=\rho_0-\rho_1$. Define
$\L_i\in H^*$ such that $\L_i(h_j)=\d_{ij},i,j=1,2.$ Then
$\rho=\L_1=-\a_2,\L_2=-\a_1-2\a_2$, and every weight can be
written as $\L=a_1\L_1+a_2\L_2$, where
$a_1=\L(h_1),a_2=\L(h_2)$ and we call $[a_1;a_2]$ the {\it Dynkin
labels} of $\L$.
\par
In this paper, we only consider (possibly infinite dimensional) $G$-modules
$V$ with finite dimensional weight space decomposition (such a module is
called a {\it weight} module): $V=\oplus_{\l\in H^*}V_\l$,
$V_\l=\{v\in V\,|\,hv=\l(h)v,\forall\,h\in H\},\,{\rm dim\ssc\,}V_\l<\infty,$
for all $\l\in H^*$.
\par
A weight vector $0\ne v_\l\in V$ with weight $\l$ is called a {\it primitive
vector} (and $\l$ is called a {\it primitive weight}) if
$v_\l\notin U(G)G_+v_\l$. If $G_+v_\l=0$, then $v_\l$ and $\l$ are called
{\it strongly primitive}. A weight $\L=a_1\L_1+a_2\L_2$ is called
{\it integral} if $a_1\in\Z$; {\it dominant} if $a_1\ge0$; {\it typical} if
$a_2\ne0\ne a_1+a_2+1$; {\it atypical} if $a_2=0$ or $a_1+a_2+1=0$ (and
accordingly, $\a_2$ or $\a_1+\a_2$ is called an {\it atypical root} of $\L$.
If we define as in [9] the atypicality matrix $A(\L)$ of a weight $\L$ to be
the $2\times1$ matrix whose entries are $A(\L)_{11}=(\L+\rho)(h_3)=
a_1+a_2+1$ and $A(\L)_{21}=(\L+\rho)(h_2)=a_2$,
then an atypical root corresponds to a zero entry in the atypicality matrix).
Note that since we shall consider possibly infinite dimensional $G$-modules,
a primitive weight $\L$ is not necessary dominant integral, therefore it
might occur that $\L$ is a multiply atypical weight (which can not occur
in finite dimensional $G$-modules). For example, when $a_1=-1,a_2=0$, then
$\L=-\L_1$ has two zero entries in the atypicality matrix. But we will call
such a weight {\it quasi-doubly atypical} instead of doubly atypical because
even such a weight occurs as a primitive weight in an indecomposable module
$V$, such a module $V$ is still not too wild in the sense of [2]. If we
investigate the Kac-module $\ol V(-\L_1)$, we see that it has three
composition factors, and it has the structure as in (2.2.i) below,
where $v_1,v_2,v_3$ are primitive vectors such that $v_1$ is the highest
weight vector of $\ol V(-\L_1),v_2=f_2v_1,v_3=f_1v_2$. However, if we use
Convention 2.3, we can regard $v_1$ to be as (2.2.ii), where $v'_1$ is the
same as $v_1$, then (2.2.i) become (2.2.iii), which fits chains' definition
in Definition 3.4.
$$
{\rm (i)}:\ \hs{10pt}
\put(0,7){$
v_1\rar v_2\put(-7,-3){$\small\vector(0,-1){8}\put(-3,-15){$v_3$}$}
$}\ \hs{100pt}\
{\rm (ii)}:\ \hs{10pt}
\put(0,7){$
v_1\put(-7,-3){$\small\vector(0,-1){8}\put(-3,-15){$v'_1$}$}
$}\ \hs{70pt}\
{\rm (iii)}:\ \hs{10pt}
\put(0,7){$
v_1\put(-7,-3){$\small\vector(0,-1){8}\put(-3,-15){$v'_1\rar$}$}
\rar v_2\put(-7,-3){$\small\vector(0,-1){8}\put(-3,-15){$v_3$}$}
$}\ \hs{70pt}\
\eqno(2.2)$$
\par
The classification of finite dimensional indecomposable weight $G$-modules
was given in Theorem 4.3 in [7], we summarize as
\par
{\bf Lemma 2.1}.
Every finite dimensional indecomposable weight $G$-module is either the
Kac-module $\ol V(\L)$ (where $\L$ is a dominant integral typical weight),
or of the form $X_4(\L),$ $X_{5a}(\L,n),$ $X_{5b}(\L,n)$ (where $\L$ is a
dominant integral atypical weight, $n\in\Z_+$), where
$X_4(\L),$ $X_{5a}(\L,n),$ $X_{5b}(\L,n)$ can be expressed by the following
diagrams (where $n$ is the number of dots):
$$
\matrix{
X_4(\L):\hfill&\cdot{\ssc\,}\matrix{\lar\cdot\lar\vs{-4pt}\cr\lar\cdot\lar\cr}
{\ssc\,}\cdot\ ,\vs{4pt}\hfill\cr
X_{5a}(\L,n):\hfill&
\cdot\rar\cdot\lar\cdots\mbox{ (ended by $\lar\cdot$ or $\rar\cdot$ if $n$ is
odd or even)},\vs{4pt}\hfill\cr
X_{5b}(\L,n):\hfill&
\cdot\lar\cdot\rar\cdots\mbox{ (ended by $\rar\cdot$ or $\lar\cdot$ if $n$ is odd
or even)}.\hfill\cr
}
\eqno\matrix{\cr(2.3)\vs{4pt}\cr\qed\cr}
$$
\par
For the purpose of next section, we also summarize below the classification
of infinite dimensional indecomposable $sl(2)$-modules given in Theorem 2.10
and Diagram 1 in [8]. For convenience, we change notations slightly: here
we denote a composition factor of a module by a dot ``$\cdot$'' just as
we did in (2.3).
\par
{\bf Lemma 2.2}.
Each indecomposable $sl(2)$-module is a subquotient module of
$U(a,b,m)$, $a,b\in\C,m\ge1$, where $U(a,b,m)$ has a basis $\{x^{(i)}_j\,|\,
j\in\Z,i\in I_m\}$, $I_m=\{1,\cdots,m\}$, with
$$
\matrix{
e_1 x^{(i)}_j=-(j+a)(j+b)x^{(i)}_{j+1}+x^{(i-1)}_{j+1},\vs{4pt}\hfill\cr
f_1 x^{(i)}_j=x^{(i)}_{j-1}+x^{(i-1)}_{j-1},\vs{4pt}\hfill\cr
h_1 x^{(i)}_j=(2j+a+b-1)x^{(i)}_j,\hfill\cr}
\eqno(2.4)$$
for $j\in\Z,i\in I_m$ (we treat $x^{(0)}_j$ as zero for $j\in\Z$).
By shifting indices of basis elements of $U(a,b,m)$ if necessary,
we can suppose $a=0,b\in\Z_+$ or $a=0,b\in\C\bs\Z$ or
$a,b\in\C\bs\Z$, and accordingly, $U(a,b,m)$ corresponds to the following
diagram:
\par
\ni\hs{15pt}
$
\put(0,0){$\put(0,8){1}\cdot$}
\put(2,0){$\small\vector(0,-1){8}$}
\put(0,-13){$\cdot$}\put(2,-13){$\small\vector(0,-1){8}$}
\put(0,-26){$\cdot$}
\put(15,-11){$\small\vector(-1,1){10}$}
\put(15,-24){$\small\vector(-1,1){10}$}
\put(15,0){$\put(0,8){2}\cdot$}
\put(17,0){$\small\vector(0,-1){8}$}
\put(15,-13){$\cdot$}\put(17,-13){$\small\vector(0,-1){8}$}
\put(15,-26){$\cdot$}
\put(15,0){$
\put(15,-11){$\small\vector(-1,1){10}$}
\put(15,-24){$\small\vector(-1,1){10}$}
$}
\put(35,-16){$\cdots$}
\put(45,0){$
\put(15,-11){$\small\vector(-1,1){10}$}
\put(15,-24){$\small\vector(-1,1){10}$}
\put(15,0){$\put(0,8){$m$}\cdot$}
\put(17,0){$\small\vector(0,-1){8}$}
\put(15,-13){$\cdot$}\put(17,-13){$\small\vector(0,-1){8}$}
\put(15,-26){$\cdot$}
$}
\put(-10,-42){(d1):\,$a=0,b\in\Z_-$}
$\ \hs{120pt}
\rb{5pt}
{$
\put(0,-13){$\put(0,8){1}\cdot$}
\put(2,-13){$\small\vector(0,-1){8}$}
\put(0,-26){$\cdot$}
\put(15,-24){$\small\vector(-1,1){10}$}
\put(15,-13){$\put(0,8){2}\cdot$}
\put(17,-13){$\small\vector(0,-1){8}$}
\put(15,-26){$\cdot$}
\put(15,0){$
\put(15,-24){$\small\vector(-1,1){10}$}
$}
\put(35,-18){$\cdots$}
\put(50,0){$
\put(15,-24){$\small\vector(-1,1){10}$}
\put(15,-13){$\put(0,8){$m$}\cdot$}
\put(17,-13){$\small\vector(0,-1){8}$}
\put(15,-26){$\cdot$}
$}
\put(-15,-45){(d2):\,$a=0,b\in\C{\ssc}\bs{\ssc}\Z$}
$}\ \hs{100pt}
\rb{-20pt}{
$\matrix{
^{\dis1}_{\dis\cdot\lar}{}^{\dis2}_{\dis\cdot\lar\,\cdots\,\cdot\lar}
{}^{\dis m}_{\dis\cdot}\vs{5pt}\cr
\mbox{(d3): $a,b\in\C{\ssc}\bs{\ssc}\Z$}\cr}
$}\hfill\rb{-15pt}{$\matrix{\hfill\cr(2.5)\cr\cr\hfill\qed\cr}$}
\par
Note that when $a=b=0$, the diagram should have the form (d2), but for
convenience of next section (mainly, Definition 3.4 and Theorem 3.5, also
{\it cf.} (2.2)$\ssc\,$), we use the following
\par
{\bf Convention 2.3}.
If $a=b=0$ in (d1), we regard all dots in the last line as the same of the
corresponding dots in the second line, i.e.,  all dots in the last line are
overlapping the dots in the second line.
\qed
\par
Also note that in (d1), $x^{(i)}_0,x^{(i)}_{-b}$ are $sl(2)$-primitive
vectors respectively with weights $b-1$ and $-b-1$, $i\in I_m$;
and in (d2), $x^{(i)}_0$ is $sl(2)$-primitive vector with weight $b-1$.
Thus an indecomposable module corresponds to a connected subdiagram
of the above diagrams (a subdiagram $W$ of diagram $U$ is a subset of dots
together with links satisfying: for any three dots $u,v,w\in U$, if
$u,w\in W$, $u$ is {\it derived by} $v$ (which means that one can go
through some way from $v$ to $u$ by following arrow's direction in diagram
$U$) and $v$ is derived by $w$, $\Rar$ $v\in W\ssc\,$).
\par\
\par
\cl{
              \S3. Infinite dimensional indecomposable modules
}
\par
Now suppose $V$ is any indecomposable $G$-module and $\l_0$ is a fixed
weight of $V$. Then $V=\oplus_{i,j\in\sZ}V_{\l_0+i\a_1+j\a_2}$. For
$j\in\Z$, let $V^{(j)}=\oplus_{i\in\sZ}V_{\l_0+i\a_1+j\a_2}$, then it is a
$G_{\ol0}$-submodule of $V$ and $e_2V^{(j)}\subset V^{(j+1)}$,
$f_2V^{(j)}\subset V^{(j-1)}$ and $V=\oplus_{j\in\sZ}V^{(j)}$ is a direct sum
of $G_{\ol0}$-submodules. Let $V^+=\{v\in V\,|\,G_{+1}v=0\}$, then $V^+$ is
a $G_{\ol0}$-submodule. By re-choosing $\l_0$ if necessary, we can suppose
$V^{+(0)}=V^+\cap V^{(0)}\ne0$. Let $W^{+(0)}$ be a maximal indecomposable
$G_{\ol0}$-submodule of $V^{+(0)}$. Since
$G_{\ol0}=\la e_1,h_1,f_1\ra\oplus\C\cc\cong sl(2)\oplus\C\cc$ is reductive
(where $\cc=h_1+2h_2$ is central in $G_{\ol0}\ssc\,$), $W^{+(0)}$ corresponds
to a subdiagram of a diagram in (2.5). Let $\wh W^{+(0)}$ be the
$G$-submodule generated by $W^{+(0)}$, then $\wh W^{+(0)}$ is a
{\it permissible quotient module} of the {\it induced module}
$\ol W{}^{+(0)}$ defined in the following.
\par
{\bf Definition 3.1}.
(1) For an indecomposable $G_{\ol 0}$-module $W$, extend it to a
$G_{\ol0}\oplus G_{+1}$-module by requiring $G_{+1}W=0$ and define
the {\it induced module} $\ol W$ of $W$ to be the $G$-module
$$
\ol W={\rm Ind\ssc\,}_{G_{\ol0}\oplus G_{+1}}^G{\ssc}W
=U(G)\otimes_{U(G_{\ol0}\oplus G_{+1})}W\cong U(G_{-1})\otimes W.
\eqno(3.1)$$
In particular, if $W=V^0(\L)$ is the highest weight $G_{\ol0}$-module
with highest weight $\L$, then the induced module $\ol W$ is the Kac-module
$\ol V(\L)$ (see, for example, [7,9,10]). Similarly, one can define the
{\it anti-induced module} $\ul{W}$ of $W$ to be the $G$-module
${\rm Ind\ssc\,}_{G_{\ol0}\oplus G_{-1}}^G{\ssc}W$ by requiring
$G_{-1}W=0$. It is straightforward to verify that $\ol W, \ul W$ are
indecomposable $G$-modules.
\par
(2) A $G$-submodule $W'$ of the induced module $\ol W$ or $\ul W$
is called a {\it permissible submodule} if $W'\cap W=0$ and accordingly,
$\ol W/W'$ or $\ul W/W'$ is called a {\it permissible quotient module}.
\par
(3) $\ol W$ or $\ul W$ is called {\it typical} if it does not contain a
nonzero permissible submodule, otherwise it is called {\it atypical}. A
subquotient module of an atypical module is also called an {\it atypical
module}.
\qed
\par
{\bf Lemma 3.2}. Let $W=W(a,b,c,m)$ be a subquotient $G_{\ol0}$-module
of $U(a,b,m)$ such that $\cc=h_1+2h_2$ acts as the scalar $c$ on $W$,
then the induced module $\ol W$ is atypical if and only if $a-b+c+1=0$ or
$a-b-c-1=0$.
\par
{\bf Proof.}
``$\Lar$'':
First suppose $W=U(a,b,m)$. If $a=0,\,b\in\Z_+$, then $x^{(1)}_{-b}$ is a
$G$-primitive vector with weight $\l$ whose Dynkin labels are
$[a_1;a_2]=[-b-1;{1\over2}(c+b+1)]$ and the $G$-submodule $U(G)v^{(1)}_{-b}$
of $\ol W$ is the (infinite dimensional) Kac-module $\ol V(\l)$ which is not
simple if and only if $\l$ is atypical, i.e., $a_2=c+b+1=0$ or $a_1+a_2+1=
c-b+1=0$. Thus if $a-b+c+1=0$ or $a-b-c-1=0$, then we can find a simple
$G$-submodule $V(\l^-)$ of $\ol V(\l)$ such that $V(\l^-)\cap W=0$, where
$\l^-$ is a weight uniquely determined by $\l$ (see Lemma 2.6 in [7]).
Thus $\ol W$ is atypical. The proof for $a=0,\,b\in\C\bs\Z$ is similar.
If $a,b\in\C\bs\Z$, it is straightforward to verify that
$\{y_j=f_2x^{(1)}_j+
{1\over2}(2j+a+b-c-1)f_3x^{(1)}_{j+1}\,|\,j\in\Z\}$
($f_3$ is defined in (2.1)$\ssc\,$)
generates a $G$-submodule having trivial intersection with $W$. In general,
by considering the induced module of a simple $G_{\ol0}$-submodule of $W$,
we can find a permissible submodule.
\par
``$\Rar$'':
As above, we only need to consider the case when $W$ is a simple sub-quotient
$G_{\ol0}$-module of $U(a,b,m)$, however in this case, it is straightforward
to check that the induced module $\ol W$ is simple if
$a{\sc\!}-{\sc\!}b{\sc\!}+{\sc\!}c{\sc\!}+{\sc\!}1{\sc\!}\ne{\sc\!}0{\sc\!}
\ne {\sc\!}a{\sc\!}-{\sc\!}b{\sc\!}-{\sc\!}c{\sc\!}-{\sc\!}1$.
\qed
\par
{\bf Proposition 3.3}.
If an indecomposable $G$-module $V$ contains a typical induced module
$\ol W$, then $V$ is an induced module.
\par
{\bf Proof.}
Let $W'$ be a maximal indecomposable $G_{\ol0}$-submodule containing $W$,
then by Lemma 3.2, $\ol W{}'$ is also typical. Then we must have
$G_{+1}W'=0$ since $W'=e_2e_3f_2f_3W'$. Thus $\ol W{}'=U(G)W'$
is a $G$-submodule of $V$. We claim that $V=\ol W{}'$. To prove this, it
remains to prove that for any $v\in V\bs\ol W{}'$, we have
$e_1v,f_1v,e_2v,f_2v\notin\ol W{}'\bs\{0\}$; this can be done case by case
and by supposing first that $W'$ has the form $U(a,b,m)$ (just as in the
proof of Lemma 3.2).
\qed
\par
Now let us go back to the discussion before Definition 3.1. By Proposition
3.3, we can suppose $\ol W{}^{+(0)}$ is atypical. We shall construct some
indecomposable $G$-modules similar to (2.3) by introducing {\it chains}.
To understand how we construct indecomposable $G$-modules, let us look the
following examples. Let $W=W(a,b,c,m)$ be as in Lemma 3.2 such that
$a-b+c+1=0$ or $a-b-c-1=0$. One can construct an anti-induced module
$\ul W'$ such that there exists an embedding $\ul\si$: $W\rar\ul W'$ with
$G_{+1}\ul\si(W)=0$ and $U(G)\ul\si(W)=\ul W'$ (such $\ul W'$ is uniquely
determined by $W$). Form the direct sum of the induced and anti-induced
modules $U=\ul W'\oplus\ol W$ and let $V$ be the $G$-submodule of $U$
generated by $\{\ul\si(w)+\ol\si(w)\,|\,w\in W\}$, where $\ol\si$:
$W\rar\ol W$ is the natural embedding (elements $\ul\si(w),\ol\si(w)$ are
called {\it top points} in $\ul W', \ol W$, and we say that $V$ is obtained
from $\ul W'$ and $\ol W$ by ``joining'' the top points of $W$). Clearly $V$
is indecomposable and we denote this module by $X_3(W)$ (similar to the
notation $X_3(\L)$ in [7]), which is just a module of the form
$X_{5b}(\L,3)$ in (2.3) if $W$ is a simple highest weight $G_{\ol0}$-module.
If we let $\wt W$ be the minimal permissible quotient module of $\ol W$,
then we can form a non-split exact sequence:
$$
0\rar \wt W\rar X_4(W)\rar X_3(W)\rar 0,
\eqno(3.2)$$
and we obtain an indecomposable $G$-module $X_4(W)$ just as we did to obtain
the module $X_4(\L)$ ({\it cf.} (2.3)$\ssc\,$). Similarly, for each such
$W$, one can construct an induced module $\ol W{}'$ and an anti-induced
module $\ul W''$ such that elements in $W$ are {\it bottom points} in
$\ol W{}'$ and $\ul W''$ (i.e., there exist embeddings
$\ol\tau:W\rar\ol W{}',\,\ul\tau:W\rar\ul W''$ such that $G_{+1}\ol\tau(W)
=0=G_{+1}\ul\tau(W)$ and $\ol\tau(W),\,\ul\tau(W)$ respectively generate
maximal permissible submodules of $\ol W{}', \ul W''\ssc\,$). Now form the
direct sum $U=\ol W{}'\oplus \ul W''$, and let $U'$ be the $G$-submodule
generated by $\{\ol\tau(w)-\ul\tau(w)\,|\,w\in W\}$, then by taking the
quotient module $V=U/U'$, we obtain an indecomposable module $V$ from
$\ol W{}'$ and $\ul W''$ by ``merging'' the bottom points of $W$. We denote
this module by $X_{5a}(W,3)$, which is similar to $X_{5a}(\L,3)$. Now we
have the following definition.
\par
{\bf Definition 3.4}.
A {\it chain} is a quadruple $(a,b,c,D)$ (when there is no confusion, we
just say $D$ is a chain) such that
\par
(i) $a,b,c\in\C$ satisfying $a-b+c+1=0$ or $a-b-c-1=0$;
\par
(ii) $D=\cup_{i\in(n_1,n_2)}D_i$, where $-\infty\le n_1\le0\le n_2\le+\infty$
and $(n_1,n_2)=\{i\in\Z\,|\,n_1<i<n_2\}$, each $D_i$ is a finite disjoint
union of nonempty subdiagrams of $U(a,b,m_i)$ for some $m_i\in\Z_+$ (if we
set $D_i=\emptyset$ for $i\notin(n_1,n_2)$, then $D=\cup_{i\in\sZ}D_i$);
\par
(iii) for each even integer $2i\in(n_1,n_2)$, there are links by
``weighted'' arrows from some elements of $D_{2i}$ (the {\it top points})
to some elements of $D_{2i\pm1}$ ({\it the bottom points}); if
$u_1,u_2\in U_{2i},v_1\in U_{2i\pm1}$ such that $u_1$ is pointed to
$v_1$ with weight $0\ne x\in\C$ (i.e.,
$u_1\rar\put(-15,6){$\sc x$}v_1$), and $u_2$ is linked to $u_1$ with
arrow pointed to $u_1$, then there exists $v_2\in D_{2i\pm1}$ such that
$u_2$ is pointed to $v_2$ with the same weight $x$ and $v_2$ is
linked to $v_1$ with arrow pointed to $v_1$, i.e., from
$\put(0,7){$
u_2\put(-7,-3){$\small\vector(0,-1){8}
\put(-3,-15){$u_1\rar\put(-13,6){$\sc x$} v_1$}$}
$}
$\ \hs{35pt}, we must have (D1):
$\put(0,7){$
u_2\put(-7,-3){$\small\vector(0,-1){8}\put(-3,-15){$u_1\rar\put(-13,6){$\sc x$}$}$}
\rar\put(-13,6){$\sc x$} v_2\put(-7,-3){$\small\vector(0,-1){8}\put(-3,-15){$v_1$}$}
$}
$\ \hs{40pt}; similarly, from
$\put(0,7){$
u_2%\put(-7,-3){$\small\vector(0,-1){8}\put(-3,-15){$u_1\rar v_1$}$}
\rar\put(-13,6){$\sc x$} v_2\put(-7,-3){$\small\vector(0,-1){8}\put(-3,-15){$v_1$}$}
$}
$\ \hs{39pt},
we must also have (D1).
The weights must also satisfy the condition that a weight $x\ne1$ can only
occur inside a circle of the {\it reduced diagram} of $D$ (where the
{\it reduced diagram} of $D$ is defined to be the diagram obtained from $D$
by replacing every part of the form (D1) by
$\cdot\rar\put(-15,6){$\sc x$}\cdot$, for example, (3.3.iii) below is
the reduced diagram of (3.3.ii)$\ssc\,$), and that there is at most one
weight $x\ne1$ inside each circle;
\par
(iv) $D$ is connected (with usual meaning).
\qed
\par
The following are examples of chains (where the link
\,$\cdot\rar\put(-14,6){$\sc\! 1$}{\ssc}\cdot\,$ has been simplified
to \,$\cdot\rar\cdot\ssc\,$):
\vs{+2pt}
\par
\ni
\hs{5pt}$
\rb{-17pt}{\rm{\small (i):}\hs{1pt}\ }
\put(0,0){
\put(3,-3){\small$\sc\vector(0,-1){8}$}
\put(3,-21){$\sc\vector(0,-1){8}$}
\put(20,-3){$\sc\vector(0,-1){8}$}
\put(20,-21){$\sc\vector(0,-1){8}$}
\put(-3,15){$\sc D_0$}\put(15,15){$\sc D'_1$}
\put(0,0){$\sc u_1\rar v'_1$}
\put(0,-18){$\sc u_2\rar v'_2$}
\put(0,-36){$\sc u_3\rar v'_3$}
}
\ \hs{32pt}\ \rb{-17pt}{\rm{\small (ii):}\hs{1pt}\ }
\put(0,0){
\put(-3,15){$\sc D_0$}
\put(0,0){$\sc u_1$}
\put(3,-3){$\vector(0,-1){8}$}
\put(0,-18){$\sc u_2$}
\put(3,-21){$\vector(0,-1){8}$}
\put(0,-36){$\sc u_3$}
\put(10,-3){$\vector(1,-1){10}$}
\put(10,-21){$\vector(1,-1){10}$}
}
\put(22,0){
\put(-3,15){$\sc D_1$}
\put(0,0){$\sc v_1$}
\put(2,-3){$\vector(0,-1){8}$}
\put(0,-18){$\sc v_2$}
\put(2,-21){$\vector(0,-1){8}$}
\put(0,-36){$\sc v_3$}
}
\put(40,0){
\put(0,-3){$\vector(-1,-3){9}$}
\put(-3,15){$\sc D_2$}
\put(0,0){$\sc w_1\rar$}
\put(3,-3){$\vector(0,-1){8}$}
\put(0,-18){$\sc w_2\rar$}
\put(3,-21){$\vector(0,-1){8}$}
\put(0,-36){${\sc w_3}$}
}
\put(60,0){
\put(-3,15){$\sc D_3$}
\put(0,0){$\sc x_1$}
\put(3,-3){$\vector(0,-1){8}$}
\put(0,-18){$\sc x_2$}
%\put(3,-21){$\vector(0,-1){8}$}
\put(0,-36){$\sc x_3\lar$}
}
\put(79,0){
\put(0,-3){$\vector(-1,-1){10}$}
\put(-3,15){$\sc D_4$}
\put(0,0){$\sc y_1\rar$}
\put(3,-3){$\vector(0,-1){8}$}
\put(0,-18){$\sc y_2$}
\put(3,-30){$\vector(0,1){8}$}
\put(0,-36){$\sc y_3\rar\put(-8,5){$\sc x$}$}
}
\put(97,0){
\put(-3,15){$\sc D_5$}
\put(0,0){$\sc z_1\lar$}
\put(3,-12){$\vector(0,1){8}$}
\put(0,-18){$\sc z_2$}
\put(3,-21){$\vector(0,-1){8}$}
\put(0,-36){$\sc z_3\lar$}
}
\put(116,0){
\put(-3,15){$\sc D_6$}
\put(0,0){$\sc z'_1\rar$}
\put(0,-36){$\sc z'_3\rar\put(-8,6){$\sc x'$}$}
}
\put(135,0){
\put(-3,15){$\sc D_7$}
\put(0,0){$\sc z''_1$}
\put(3,-12){$\vector(0,1){8}$}
\put(0,-18){$\sc z''_2$}
\put(3,-21){$\vector(0,-1){8}$}
\put(0,-36){$\sc z''_3$}
}
\ \hs{155pt}\rb{-17pt}{\rm{\small (iii):}\hs{1pt}\ }
\put(0,-7)
{$
%%%%%%%%%%%%%%
\put(0,0){
\put(0,0){$\sc\cdot\rar$}
\put(2,0){$\vector(0,-1){8}$}
\put(0,-14){$\sc\cdot$}
}
\put(13,0){
\put(0,0){$\sc\cdot\lar$}
\put(2,13){$\vector(0,-1){8}$}
\put(0,13){$\sc\cdot$}
}
\put(26,0){
\put(0,0){$\sc\cdot\rar$}
\put(2,0){$\vector(0,-1){8}$}
\put(0,-14){$\sc\cdot$}
\put(2,-14){$\vector(0,-1){8}$}
\put(0,-27){$\sc\cdot$}
}
\put(39,0){
\put(0,0){$\sc\cdot\rar$}
\put(0,-27){$\sc\cdot\lar$}
}
\put(52,0){
\put(0,0){$\sc\cdot\rar$}
\put(2,0){$\vector(0,-1){8}$}
\put(0,-14){$\sc\cdot$}
\put(2,-22){$\vector(0,1){8}$}
\put(0,-27){$\sc\cdot\rar\put(-6,5){$\sc x$}$}
}
\put(65,0){
\put(0,0){$\sc\cdot\lar$}
\put(2,-9){$\vector(0,1){8}$}
\put(0,-14){$\sc\cdot$}
\put(2,-14){$\vector(0,-1){8}$}
\put(0,-27){$\sc\cdot\lar$}
}
\put(78,0){
\put(0,0){$\sc\cdot\rar$}
\put(0,-27){$\sc\cdot\rar\put(-9,5){$\sc x'$}$}
}
\put(90,0){
\put(0,0){$\sc\cdot$}
\put(2,-9){$\vector(0,1){8}$}
\put(0,-14){$\sc\cdot$}
\put(2,-14){$\vector(0,-1){8}$}
\put(0,-27){$\sc\cdot$}
}
%%%%%%%%%%%%%%%%
$}
\hfill\rb{-17pt}{(3.3)}$\vs{2pt}\par
Now we can obtain our main result in the following.
\par
{\bf Theorem 3.5}.
(1) To each chain $(a,b,c,D)$, there exists a unique indecomposable
$G$-module $X(a,b,c,D)$ corresponding to it.
(2) An indecomposable $G$-module $V$ is either a typical induced module
$\ol W$, or of the form $X_4(W),\,X(a,b,c,D)$.
\par
{\bf Proof.} (1).
For each chain $(a,b,c,D)$, to obtain $X(a,b,c,D)$, we shall use the method
above for constructing $X_3(W)$ and $X_{5a}(W,3)$. To see the pattern, take
(3.3.ii) as the example. For $i\in\Z$, let
$D^{\pm}_{2i-1}=\{v\in D_{2i-1}\,|\,\exists{\ssc\,}
u\in D_{2i-1\pm1}$, $u$ is pointed to $v\}$.
If $D$ is as in (3.3.ii), then $D^-_1=\{v_2,v_3\},$ $D^+_1=\{v_3\}$.
Below we shall first construct $G$-modules $U,V',W$ which have diagrams
respectively corresponding to $D_0\cup D^-_1,D_1,D^+_1\cup D_2$, then we
will ``merge'' all {\it bottom points} to obtain a $G$-module $M'$ whose
diagram corresponds to $D_0\cup D_1\cup D_2$.
\par
First, let $U^0=U^0(a,b,c,m_0)$ be the $G_{\ol0}$-module corresponding
to the diagram $D_0$ (if $D_0$ is a disjoint union of more than one
diagrams, then $U^0$ is not indecomposable). Let $\ol U{}^0$ be the induced
module of $U^0$ (then $\ol U{}^0$ has a picture (3.3.i) if $D$ is
(3.3.ii)$\ssc\,$). Let $\wt U{}^0$ be the maximal permissible $G$-submodule
of $\ol U{}^0$. So $\wt U{}^0$ has a picture $D'_1$ in the example, which is
in fact the same as $D_0$. Then $D^-_1$ corresponds to a quotient $G$-module
of $\wt U{}^0$ (by removing $v'_3$ in the example). Now take $U$ to be the
quotient of $\ol U{}^0$ by removing $v'_3$ so $U$ has the diagram
$D_0\cup D^-_1=D_0\cup D_1\bs\{v_1\}$ in (3.3.ii).
\par
Next, take $U^1=U^1(a_1,b_1,c_1,m_1)$ to be the $G_{\ol0}$-module
corresponding to the diagram $D_1$, where $a_1,b_1,c_1$ can be uniquely
determined by the diagram of $U$ (in above example, $a_1,b_1,c_1$ for
$U^1$ must be the same as those $a,b,c$ corresponding to the
$G_{\ol0}$-module $D_1\bs\{v_1\}$). Then by taking $V'$ to be the minimal
permissible quotient module of the induced module $\ol U{}^1$, we can obtain
a $G$-module $V'$ which has diagram $D_1$. Similarly, one can uniquely
construct a $G$-module $W$ whose diagram is
$D^+_1\cup D_2=D_1\cup D_2\bs\{v_1,v_2\}$ in the example. Now form the
direct sum $M=U\oplus V'\oplus W$ and take $M'$ to be the quotient module by
``merging'' all the bottom points (i.e., if we set elements $v_2,v_3$ in $U$
to be $v'_2,v'_3$ and element $v_3$ in $W$ to be $v''_3$, and let $M''$ be
the $G$-submodule of $M$ generated by elements $v'_2-v_2,v'_3-v_3,v_3-v''_3$
(if, say, the link from $u_2$ to $v_2$ has a weight $x{\ssc}\ne{\ssc}1$,
then $v'_2-v_2$ should be replaced by $v'_2-xv_2$), then
$M'=M/M''$), then $M'$ has diagram $D_0\cup D_1\cup D_2$. Using this way,
we can construct indecomposable module corresponding to the chain
$(a,b,c,D)$.
\par
(2) If $V$ contains $X_4(W')$ for some indecomposable $G_{\ol0}$-module
$W'$, then as in the proof of Proposition 3.3, $V$ must be $X_4(W)$, where
$W$ is the maximal $G_{\ol0}$-submodule containing $W'$. Now suppose $V$
is not a typical induced module, nor is of the form $X_4(W)$. Let $V'$ be
a maximal $G$-submodule of $V$ which has the form $X(a,b,c,D)$ for some
chain $(a,b,c,D)$, then 
one can prove that $V=V'=X(a,b,c,D)$.
(Take $D$ to be a suitable minimal generating set of $V$ so that each
element of $D$ corresponds to a composition factor of $V$,
defining a relation on $D$ by: for $u\ne v\in D,
u\rar v\Lra v\in U(G_-)G_{-1}u$;
$u\lar v\Lra u\in U(G_+)G_{+1}v$;
$\put(0,7){$u$}\put(3,4){$\vector(0,-1){7}$}\put(0,-9){$v$}$\hs{6pt}
$\Lra v\in U(G_-)f_1u$ and
if $w\in D,w\in U(G_-)f_1u,
v\in U(G_-)f_1w$, then $v=w$;
$\put(0,7){$u$}\put(3,-3){$\vector(0,1){7}$}\put(0,-9){$v$}$\hs{6pt}
$\Lra u\in U(G_+)e_1v$ and if $w\in D,w\in U(G_+)e_1v,
u\in U(G_+)e_1w$, then $u=w$; then one can prove that by a suitable
choice\,of $D$, it is a chain.)
\qed
\par\
\par
\cl{                     ACKNOWLEDGEMENTS
}
The author was supported by a fund from National Education Department of
China.
\par\
\par
\cl{                      REFERENCES
}
\par
\ni
\baselineskip 6pt
\lineskip 6pt
\hi3.5ex\ha1
[1] Chmelev G.~C. Indecomposable representations of $sl(1,2)$.
 Comptes Rendus de l'Acad\'emie Bulgare des Sciences, 
{\bf 1982}, {\it 35}(8), (in russian).
\par
\ni
\hi3.5ex\ha1
[2] Germoni J. Indecomposable representations of special linear Lie
superalgebras. J.~Alg. {\bf 1998}, {\it 209}, 367--401.
\par
\ni
\hi3.5ex\ha1
[3] Kac V.~G. Lie superalgebras. Adv.~in Math. {\bf1977}, {\it26}, 8--96.
\par
\ni
\hi3.5ex\ha1
[4] Kac V.~G. Representations of classical Lie superalgebras.
Lect.~in Math., Springer-Verlag, Berlin, ed.~Bleuler K., Petry H.~and
Reetz A., {\bf 1977}, {\it 676}, 579--626.
\par
\ni
\hi3.5ex\ha1
[5] Leites D. Indecomposable representations of $sl(1|n)$.
Theoretical Physics and Mathematics, Acad. Nauk. CCCP, {\bf 1982},
(in russian).
\par
\ni
\hi3.5ex\ha1
[6] Su Y. A complete classification of indecomposable Harish-Chandra
modules of the Lie superalgebra $B(0,1)$.
Comm.~Alg. {\bf1992}, {\it20}, 573--582.
\par
\ni
\hi3.5ex\ha1
[7] Su Y. Classification of finite dimensional modules of the Lie
superalgebra $sl(2/1)$. Comm.~Alg. {\bf1992}, {\it20}, 3259--3277.
\par
\ni
\hi3.5ex\ha1
[8] Su Y. A classification of indecomposable $sl_2(\C)$-modules
 and a conjecture of Kac on Irreducible modules over the Virasoro
 algebra. J.~Alg. {\bf1993}, {\it161}, 33--46.
\par
\ni
\hi3.5ex\ha1
[9] Su Y.; Hughes J.~W.~B.; King R.~C. Primitive vectors of Kac-modules
of the superalgebras $sl(m/n)$.
J.~Math.~phys. {\bf2000}, {\it41}, 5064--5087.
\end{document}